\documentclass[12pt]{article}
\usepackage{amsmath,amssymb}
\usepackage{graphicx,psfrag,epsf}
\usepackage{enumerate}
\usepackage[sort, numbers]{natbib}
\usepackage{url} 
\usepackage{multicol}
\usepackage{subfig}
\usepackage{url}

\graphicspath{{./}{./figs/}}

  \def\clap#1{\hbox to 0pt{\hss#1\hss}}

\usepackage{bm}
\providecommand{\mat}[1]{\bm{#1}}%
\renewcommand{\vec}[1]{\mathbf{#1}}

\providecommand{\mC}{\ensuremath{\mat{C}}}

\providecommand{\mW}{\ensuremath{\mat{W}}}

\providecommand{\ve}{\ensuremath{\vec{e}}}

\providecommand{\vw}{\ensuremath{\vec{w}}}
\providecommand{\vx}{\ensuremath{\vec{x}}}
\providecommand{\vy}{\ensuremath{\vec{y}}}
\providecommand{\vz}{\ensuremath{\vec{z}}}

\newcommand{\hmC}{\hat{\mC}}

\newcommand{\hLambda}{\hat{\Lambda}}

\newcommand{\hmW}{\hat{\mW}}
\newcommand{\hvw}{\hat{\vw}}

\newcommand{\sL}{\mathcal{L}}

\newcommand{\sU}{\mathcal{U}}
\newcommand{\sX}{\mathcal{X}}

\newcommand{\bpi}{\bm{\pi}}

\newcommand{\bmat}[1]{\begin{bmatrix}#1\end{bmatrix}}

\newcommand{\maxi}[1]{\underset{#1}{\mathrm{maximum}}}

\newcommand{\IL}{I_{\mathrm{L}}}
\newcommand{\IS}{I_{\mathrm{S}}}
\newcommand{\RS}{R_{\mathrm{S}}}
\newcommand{\RP}{R_{\mathrm{P}}}
\newcommand{\ISC}{I_{\mathrm{SC}}}
\newcommand{\NS}{N_{\mathrm{S}}}
\newcommand{\Vth}{V_{\mathrm{th}}}

\newcommand{\Pmax}{P_{\text{max}}}

\newtheorem{theorem}{Theorem}[section]
\newtheorem{lemma}[theorem]{Lemma}

\begin{document}

\def\spacingset#1{\renewcommand{\baselinestretch}%
{#1}\small\normalsize} \spacingset{1}


  \title{\bf Discovering an Active Subspace in a\\ Single-Diode Solar Cell Model}
  \author{Paul G.~Constantine\thanks{
    The authors gratefully acknowledge the support Department of Energy's Advanced Scientific Computing Research Award \emph{Mathematical and Statistical Methodologies for DOE Data-Centric Science at Scale}.}\hspace{.2cm}
    and 
    Brian Zaharatos\\
    Department of Applied Mathematics and Statistics\\ 
    Colorado School of Mines\\
    and \\
    Mark Campanelli\\
    National Renewable Energy Laboratory}
  \maketitle

\bigskip
\begin{abstract}
Predictions from science and engineering models depend on the values of the model's input parameters. As the number of parameters increases, algorithmic parameter studies like optimization or uncertainty quantification require many more model evaluations. One way to combat this curse of dimensionality is to seek an alternative parameterization with fewer variables that produces comparable predictions. The \emph{active subspace} is a low-dimensional linear subspace defined by important directions in the model's input space; input perturbations along these directions change the model's prediction more, on average, than perturbations orthogonal to the important directions. We describe a method for checking if a model admits an exploitable active subspace, and we apply this method to a single-diode solar cell model with five input parameters. We find that the maximum power of the solar cell has a dominant one-dimensional active subspace, which enables us to perform thorough parameter studies in one dimension instead of five. 
\end{abstract}

\noindent%
\vfill
\hfill {\tiny CoDA 2014 {\em Statistical Analysis and Data Mining} tex template (do not remove)}

\newpage
\spacingset{1.45} 

\section{Introduction}
\label{sec:intro}

Science and engineering simulations often contain several input parameters---e.g., physical constants, boundary conditions, or geometry descriptions. When presented with such parameterized simulations, the scientist naturally wonders how the simulated predictions depend on the input parameters. Which parameters, when perturbed, create the largest change in predictions? How precisely must the parameters be specified to ensure accurate predictions? And what is the effect of imprecisely prescribed input parameters on the predictions? If the goal is to maximize or minimize a predicted quantity, which combinations of parameters correspond to larger or smaller values of the prediction? In a similar vein, one may ask which parameter values correspond to predictions outside a safe region of operation---or which parameter values yield predictions that are consistent with a set of observations. 

Many scientists rely on intuition about the physical system to answer these questions. But intuition becomes less trustworthy as the simulations become more complex, e.g., when they include several interacting physical components. Algorithms for optimization, uncertainty quantification, and model calibration become more attractive as model complexity increases. If an algorithm can easily interface with the simulation---e.g., by automatically evaluating predictions given values for the inputs---then applying the algorithm to the simulation becomes relatively easy. 

The number of times the algorithm needs to evaluate a prediction increases---sometimes extremely rapidly---as the number of inputs increases. The situation is worse if each evaluation requires significant computational resources. For example, finding the global minimum of a complicated prediction depending on 100 inputs is not tractable if available resources permit only ten model evaluations. In practice, the scientist may choose only the most important parameters to vary so the study fits within the computational budget. Alternatively, one may seek a low-dimensional description of the prediction as a function of the input parameters. If such a description is sufficiently accurate, then studies like optimization or calibration can work in the space of fewer variables---potentially allowing the desired parameter studies within the given budget. 

The \emph{active subspace} is a low-dimensional linear subspace defined by important directions in the model's input space; input perturbations along these directions change the model's prediction more, on average, than perturbations orthogonal to the important directions. Not all models have an active subspace. Some model's predictions respond significantly to input perturbations along all directions. However, if a model does admit an active subspace, then one can exploit it to perform parameter studies in the coordinates of the subspace---i.e., the \emph{active variables}---which are linear combinations of normalized versions of the model's input parameters. Therefore, it can be very advantageous to discover that a model admits an active subspace. 

Active subspaces have been studied in a variety of contexts under different names. Cook's excellent text \emph{Regression Graphics}~\cite{cook2009regression} reviews and develops statistical methods for dimension reduction in the context of regression surfaces, and it contains references to the major works in the statistics literature. What we call the \emph{active subspace} is a type of \emph{dimension-reduction subspace} in Cook's parlance~\cite[Chapter 6]{cook2009regression}---though we are working with noiseless computer simulations in contrast to general regression surfaces. Russi's 2010 Ph.D.~dissertation uses the phrase \emph{active subspace} in a way comparable to our use~\cite{Russi2010}. He exploits the active subspace to construct quadratic surrogate models for uncertainty quantification in chemical kinetics models. Our prior work develops a theoretical framework for active subspaces including applications to kriging response surfaces~\cite{Constantine2013}. We have applied these methods to several models in aerospace engineering~\cite{constantine2011,wang2011,lukaczyk2014}. Abdel-Khalik has applied similar methods in nuclear engineering applications~\cite{AbdelKhalik2013}.

In this paper, we describe how to test a model for an active subspace, and we apply this test to the maximum power from a single-diode model of a photovoltaic solar cell with five input parameters. In Section \ref{sec:asm}, we generically describe the active subspace and how to search for it. We then describe the single-diode model, its input parameters, and its predicted \emph{performance parameters} (i.e., model outputs) in Section \ref{sec:pv}. We apply the tests for the active subspace to the model's maximum power in Section \ref{sec:pvasm} and show that a dominant one-dimensional active subspace is present in the five-dimensional space of input parameters. We conclude in Sections \ref{sec:discussion} and \ref{sec:conclusion} with a discussion of how one may exploit the low-dimension of the active subpace to further study the parameter dependence in the single-diode model's maximum power. 

\section{Active subspaces}
\label{sec:asm}


We consider a generic multivariate function $f=f(\vx)$, where $\vx$ represents the inputs of the model, and $f$ represents a specific scalar performance parameter that the model predicts. Let $\sX=[-1,1]^m$ be the domain with $\vx\in\sX$, and let $\rho:\sX\rightarrow\mathbb{R}_+$ be a bounded and continuous weight function on $\sX$; we assume $\rho$ is normalized to integrate to 1. Assume $f$ is diffentiable and absolutely continuous, and denote the gradient $\nabla_\vx f(\vx) = [\partial f/\partial x_1,\dots,\partial f/\partial x_m ]^T$ oriented as a column vector. 

Consider the following matrix $\mC$ defined as
\begin{equation}
\label{eq:C}
\mC \;=\; \int_\sX (\nabla_\vx f)(\nabla_\vx f)^T\,\rho\,d\vx.
\end{equation}
In the context of dimension reduction for regression functions, Samarov calls this matrix an \emph{average derivative functional}~\cite{samarov1993exploring}. Note that we are not studying regression functions, per se. Instead, $f$ represents a noiseless, parameterized computer simulation. The matrix $\mC$ is symmetric and positive semi-definite, so it admits a real eigenvalue decomposition
\begin{equation}
\mC \;=\; \mW\Lambda\mW^T, 
\qquad 
\Lambda \;=\; \mathrm{diag}(\lambda_1,\dots,\lambda_m),
\qquad
\lambda_1\geq\cdots \geq\lambda_m\geq 0.
\end{equation}
We can partition the eigenvectors,
\begin{equation}
\label{eq:partition}
\mW = \bmat{\mW_1 & \mW_2},\qquad \Lambda = \bmat{\Lambda_1 & \\ & \Lambda_2},
\end{equation}
where $\mW_1$ contains the first $n$ eigenvectors, and $\Lambda_1$ contains the $n$ largest eigenvalues. We use the two sets of eigenvectors to create new sets of variables $\vy=\mW_1^T\vx$ and $\vz=\mW_2^T\vx$. We call the subspace defined by $\mW_1$ the \emph{active subspace}, and we call the variables $\vy$ the \emph{active variables}; the term \emph{active subspace methods} for this type of analysis was first used in Russi's 2010 Ph.D.~dissertation~\cite{Russi2010}. The following two lemmas justify these labels.

\begin{lemma}
\label{lem:avgsqgrad}
The mean-squared directional derivative of $f$ with respect to the eigenvector $\vw_i$ is equal to the corresponding eigenvalue,
\begin{equation} 
\int_\sX \big((\nabla_{\vx} f)^T\vw_i\big)^2\,\rho\,d\vx \;=\; \lambda_i.
\end{equation}
\end{lemma}

\begin{lemma}
\label{lem:grad}
The mean-squared gradients of $f$ with respect to the coordinates $\vy$ and $\vz$ satisfy
\begin{equation}
\begin{aligned}
\int_\sX (\nabla_{\vy} f)^T(\nabla_{\vy} f)\,\rho\,d\vx &= \lambda_1+\cdots+\lambda_n, \\
\int_\sX (\nabla_{\vz} f)^T(\nabla_{\vz} f)\,\rho\,d\vx &= \lambda_{n+1}+\cdots+\lambda_m.
\end{aligned}
\end{equation}
\end{lemma}
The proofs of these lemmas can be found in our prior work~\cite{Constantine2013}. In words, they mean that $f$ changes more, on average, when its inputs are perturbed along the directions $\mW_1$ than along the directions $\mW_2$; the eigenvalues quantify precisely how much more. If an eigenvalue is exactly zero, then $f(\vx)$ is constant along the direction defined by the corresponding eigenvector over all of $\sX$. 

To gain some intuition, consider the extreme case where all eigenvalues are precisely zero except $\lambda_1$. Then $f(\vx) = g(\vw_1^T\vx)$, where $\vw_1$ is the first column of $\mW$, and $g$ is a function of one variable. In many applications, the smaller eigenvalues are not precisely zero, but they may be much (e.g., orders of magnitude) smaller so that $f(\vx)$ may be reasonably approximated by a function of $n<m$ linear combinations of $\vx$. 

If a given model $f$ admits such structure, then certain operations---e.g., response surface modeling or optimization---become much less expensive. In particular, these operations can be performed in the $n$-dimensional space of the active variables $\vy$ instead of the full $m$-dimensional space. It is therefore extremely valuable to determine if $f$ admits an active subspace. This analysis assumes that the weight function $\rho$ on the domain $\sX$ is given. The computed quantities (like $\mW$ and $\Lambda$) will change if a different $\rho$ is given. 

\subsection{Identifying an active subspace}
To identify the active subspace, we must approximate the matrix $\mC$ from \eqref{eq:C}. Since $\mC$ is the mean of the outer product of the gradient, we can approximate it with simple Monte Carlo. One could achieve a more accurate approximation of $\mC$ and its eigendecomposition with an integration rule that is more accurate than simple Monte Carlo. However, if $m$ is greater than two or three, then tensor product constructions of accurate univariate numerical integration rules (e.g., Gaussian quadrature) require too many evaluations of the gradient to be practical---especially if the gradient is expensive to compute. The simple Monte Carlo proceeds as follows. Draw $\vx_i$ independently according to $\rho$ with $i=1,\dots,M$. Then
\begin{equation}
\label{eq:Cmc}
\mC \;\approx\;\hmC\;=\;
\frac{1}{M}\sum_{i=1}^M \nabla_\vx f(\vx_i)\,\nabla_\vx f(\vx_i)^T \;=\;
\hmW\hLambda\hmW^T.
\end{equation}
We use the eigenvalues of the approximation $\hmC$ as evidence of an active subspace. In particular, a large gap in the eigenvalues indicates a separation between the corresponding active and inactive subspaces defined by $\hmW_1$ and $\hmW_2$, respectively. The qualification \emph{large} depends on the application. For example, a good low-dimensional approximation of $f$ over its entire domain may need a much larger separation than a good approximation of the bounds or the average of $f$. 

Other eigenvalue-based dimension reduction schemes (e.g., principal component analysis) use heuristics based on the magnitude of the eigenvalues, such as choosing the dimension of the subspace such that the ratio $\sum_{i=1}^n \lambda_i/\sum_{i=1}^m \lambda_i$ is larger than some threshold. For active subspaces, such heuristics are not well-justified. We are not interested in capturing the gradient's variance; we want to identify directions along which the function changes the most. In our recent work, we show that the distance between the $n$-dimensional subspace defined by $\mW_1$ and its approximation defined by $\hmW_1$ is inversely proportional to $\lambda_{n+1}-\lambda_n$~\cite{constantine2014computing}. Therefore, the quality of the finite-sample subspace approximation depends more on the gap between the eigenvalues than their magnitude. 

\subsection{Bootstrap to estimate variability}
\label{sec:boot}
If a gap is present in the approximate eigenvalues $\hLambda$, one may naturally ask if a comparable gap is present in the true eigenvalues $\Lambda$. To address this question, we use nonparametric bootstrap. The bootstrap is most appropriate when the simulation is expensive, and the number $M$ of gradient evaluations is constrained by a computational budget. For $j=1,\dots,M'$, let $\bpi_j = [\pi_1^j,\dots,\pi_M^j]$ be an $M$-vector of integers drawn uniformly at random between 1 and $M$. The $j$th \emph{bootstrap replicate} of $\hmC$ is computed as
\begin{equation}
\label{eq:Cboot}
\hmC^j\;=\;
\frac{1}{M}\sum_{i=1}^M \nabla_\vx f(\vx_{\pi_i^j})\,\nabla_\vx f(\vx_{\pi_i^j})^T \;=\;
\hmW^j\hLambda^j\left(\hmW^j\right)^T.
\end{equation}
The collection of eigenvectors $\hLambda_j$ yields a bootstrap distribution for the eigenvalues, which can be used to estimate boostrap intervals. Note that these intervals are not true confidence intervals, since the estimates are biased. However, this bias decreases as $M$ increases; see~\cite[Section 7.2]{efron1994introduction} and~\cite[Chapter 3]{jolliffe2002principal} for related discussions of the bias in bootstrap estimates for principal components. Nevertheless, we use the 99\% bootstrap intervals as evidence of gaps in the estimated eigenvalues. 

We can also use the bootstrap replicates to estimate the error in the estimated subspace. Partition the eigenvectors from the bootstrap replicate as in \eqref{eq:partition}. Define
\begin{equation}
\label{eq:subspacerr}
e_j \;=\; \|\hmW_1\hmW_1^T - (\hmW^j_1)(\hmW^j_{1})^T\|_2
\;=\; \|\hmW_1^T\hmW^j_{2}\|,
\end{equation}
which is the distance between the subspace defined by $\hmW_1$ and bootstrap replicate $\hmW_1^j$~\cite{golub1996matrix}. The mean and 99\% bootstrap intervals of the set $\{e_j\}$ quantify the varibility in the estimated active subspace. 

Our current research efforts are focused on rigorously justified, computable metrics for determining the relationship between the approximate eigenpairs $\hLambda$, $\hmW$ and the true eigenpairs $\Lambda$, $\mW$. Our recent paper applies non-asymptotic random matrix theory to characterize the error in these approximations~\cite{constantine2014computing}.

\subsection{Visualization with the active subspace}
Scatter plots are a common way to visualize data sets in search of a trend. Unfortunately, visualization tools can only display scalar responses $f$ as a function of at most two variables. When the response depends on more than two variables, one can plot responses versus each variable or each pair of variables. 

Using the vectors defining the active subspace, we can create scatter plots based on the active variables, which often show discernible trends that can be exploited in further studies (e.g., building response surfaces or optimization). This idea is described in detail in the texts on \emph{regression graphics}~\cite{cook2009regression}, where such plots are called \emph{sufficient summary plots}. For example, suppose we have noticed a large gap between the first and second eigenvalue. We can confirm the presence of the active subspace by first sampling $\vx_i$ from $\rho$ as in \eqref{eq:Cmc}, computing $f_i = f(\vx_i)$ (often computed along with the gradients in \eqref{eq:Cmc}), and plotting the pairs $(\vw_1^T\vx_i,f_i)$---where $\vw_1$ is the first column of $\mW$. If a tight, univariate trend is visibly present, then this verifies the active subspace. Figure \ref{fig:proj1} in Section \ref{sec:pvasm} shows these scatter plots for the single-diode model.

These sufficient summary plots can also be used to examine the variability in the vectors defining the active subspace. In particular, we can use the bootstrap replicates $\hmW^j$ from \eqref{eq:Cboot}. The first eigenvector from each can be used to create pairs for the scatter plot. Plotting all pairs together creates clusters in the scatter plot that can indicate the spread around a trend. We create these visualizations for the predicted performance parameter from the single-diode model in Section \ref{sec:pvasm}.

\subsection{Gradient approximation}
Often the gradient $\nabla_\vx f$ is not available or is too complicated to compute. To approximate partial derivatives, one builds a model of $f(\vx)$ that is easily differentiable, e.g., a polynomial model. A first-order finite difference approximation of the partial derivative with respect to $x_i$ at $\vx$ computes the slope of the plane that interpolates $f$ at $\vx$ and $\vx+\varepsilon \ve_i$, where $\ve_i$ is a vector of zeros with a one in the $i$th component. When the variables are properly normalized, a reasonable choice for $\varepsilon$ to approximate the derivative and avoid round-off issues is the square root of machine precision~\cite[Section 8.1]{Nocedal2006}.

If the simulation output is well-behaved as a function of the variables $\vx$, then finite difference approximations of gradients can be used as a substitute for the true gradient. Each finite difference approximation requires $m+1$ evaluations of $f$, so the number of evaluations to compute $\hmC$ is $M(m+1)$. Therefore, finite differences are most appropriate when the following are satisfied:
\begin{itemize}
\item the simulation output behaves well as a function of the parameters, i.e., there is no noise due to limited iterations of a nonlinear solver or large changes of the function on scales smaller than the finite difference step size,
\item the simulation completes in a short enough time on available computing resources to permit $M(m+1)$ evaluations. 
\end{itemize}
Both of these conditions are satisfied for the performance parameter of interest from the single-diode model as a function of its input parameters. Therefore, we use finite differences in place of the true gradient to search for an active subspace. 

If these conditions are not satisfied, then one must pursue other options for estimating the gradients. One idea is to use local linear models within a budget-constrained set of runs. In short, assume one has $B$ pairs $(\vx_j,f(\vx_j))$. Given a point $\vx$ in the parameter space, select $M$ points in $\{\vx_j\}$ near $\vx$, fit a linear model to the selected pairs, and estimate the gradient as the gradient of the linear model. $M$ must be large enough and the points must be chosen to enable fitting the linear model. We have used this approach in practice, but we are still analyzing its properties. The advantage over finite differences is that it is robust to noise in $\{f_j\}$ when the number of points used to fit the linear model is greater than the number of variables ($M>m+1$). However, if the points $\{\vx_j\}$ are not sufficiently dense in the parameter space, then the gradient estimates can be very poor, which lead to poor estimates of the eigenvalues and eigenvectors. We continue to analyze this approach. 

\section{A single-diode solar cell model}
\label{sec:pv}

In this section we introduce the single-diode solar cell model, the input parameters, and the predicted quantity of interest. Over the past decade, a growing demand for clean energy has caused rapid growth in the photovoltaic (PV) industry. The overall health of the industry depends on proper characterization of the risks associated with a PV system. One can characterize such risks using a mathematical model to describe the performance of the PV device, such as a single-diode lumped-parameter equivalent-circuit model defined in Section \ref{sec:params}. For many series-wired PV devices, the single-diode model accurately describes the device's current-voltage ($I$-$V$) characteristics at given irradiance and temperature conditions. The single-diode model contains several input parameters that must be estimated---e.g., with methods in \cite{Hansen2013b} or \cite{Farivar2010}---to properly pose the model. With the inputs fixed, practitioners can estimate the {\it key performance parameters}, such as the maximum power output $\Pmax$ and energy conversion efficiency $\eta$. 

Estimating the model's input parameters presents challenges. For example, many single-diode model parameter estimation methods have been proposed in the PV literature, and it is not always clear if any particular estimation yields a unique result \cite{Hansen2013}. Further, since each key performance parameter is a function of the model parameters, errors in estimating the latter will cause issues in estimating the former. 

Given these issues, it is desirable to know how the variability or uncertainty in each of the input parameters affects variability in the outputs. In particular, it would be useful to know whether the model admits an active subspace. If an active subspace is present, then estimation of the performance parameters is greatly simplified. For example, rather than conducting current-voltage measurements in a way that yields accurate estimates of all model parameters, practitioners could instead focus on measurements that accurately estimate parameters in the active subspace. 

\subsection{Calculating the key performance parameters}
\label{sec:params}

The single-diode model describes a relationship between the current ($I$) and voltage ($V$) in the single diode. There are multiple variants of the single-diode model, each of which has accompanying auxiliary equations that describe how components depend on irradiance and temperature. In this work, we use the auxiliary equations derived in \cite{Zaharatos2014}; see \cite{DeSoto2006} for an example of others. These relationships are defined in the following equations.
\begin{equation}
\label{eq:model}
I  \;=\; \IL - \IS\,\left(\exp\left(\frac{V + I\,\RS}{\NS\,n \,\Vth}\right) - 1 \right) - \frac{V + I\,\RS}{\RP}. 
\end{equation}
$\NS$ is the number of cells connected in series, which we set to 1. The thermal voltage $\Vth$ is a given, fixed constant for fixed temperature $T=25^\circ\mathrm{C}$. An auxiliary equation defines the photocurrent $I_\text{L}$ as
\begin{equation}
\label{eq:photocurrent}
I_\text{L} \;=\; \ISC + \IS \left(\exp\left(\frac{\ISC \RS}{\NS\, n\, \Vth}\right) - 1 \right) + \frac{\ISC\,\RS}{\RP}.
\end{equation}
The remaining terms are input parameters to the single-diode model; their names, units, and ranges are given in Table \ref{tab:domains}. The ranges correspond to typical 2cm$^2$ crystalline silicon PV cells.  

\begin{table}[ht]\small
\caption{The ranges of the input parameters for the single-diode model in \eqref{eq:model} and \eqref{eq:photocurrent}. These ranges correspond to typical 2cm$^2$ crystalline silicon PV cells.}
\begin{center}
\begin{tabular}{lllll}
Parameter & Name & Lower bound & Upper bound & Units\\
\hline
$\ISC$ & short-circuit current & 0.05989 & 0.23958 & amps\\
$\IS$ & diode reverse saturation current & 2.2e-11 & 2.2e-7 & amps \\
$n$ & ideality factor & 1 & 2 & unit-less\\
$\RS$ & series resistance & 0.16625 &  0.66500 & ohms \\
$\RP$ & parallel (shunt) resistance & 93.75 & 375.00 & ohms
\end{tabular}
\end{center}
\label{tab:domains}
\end{table}

\begin{figure}[ht]
\begin{center}
\includegraphics[width=0.45\textwidth]{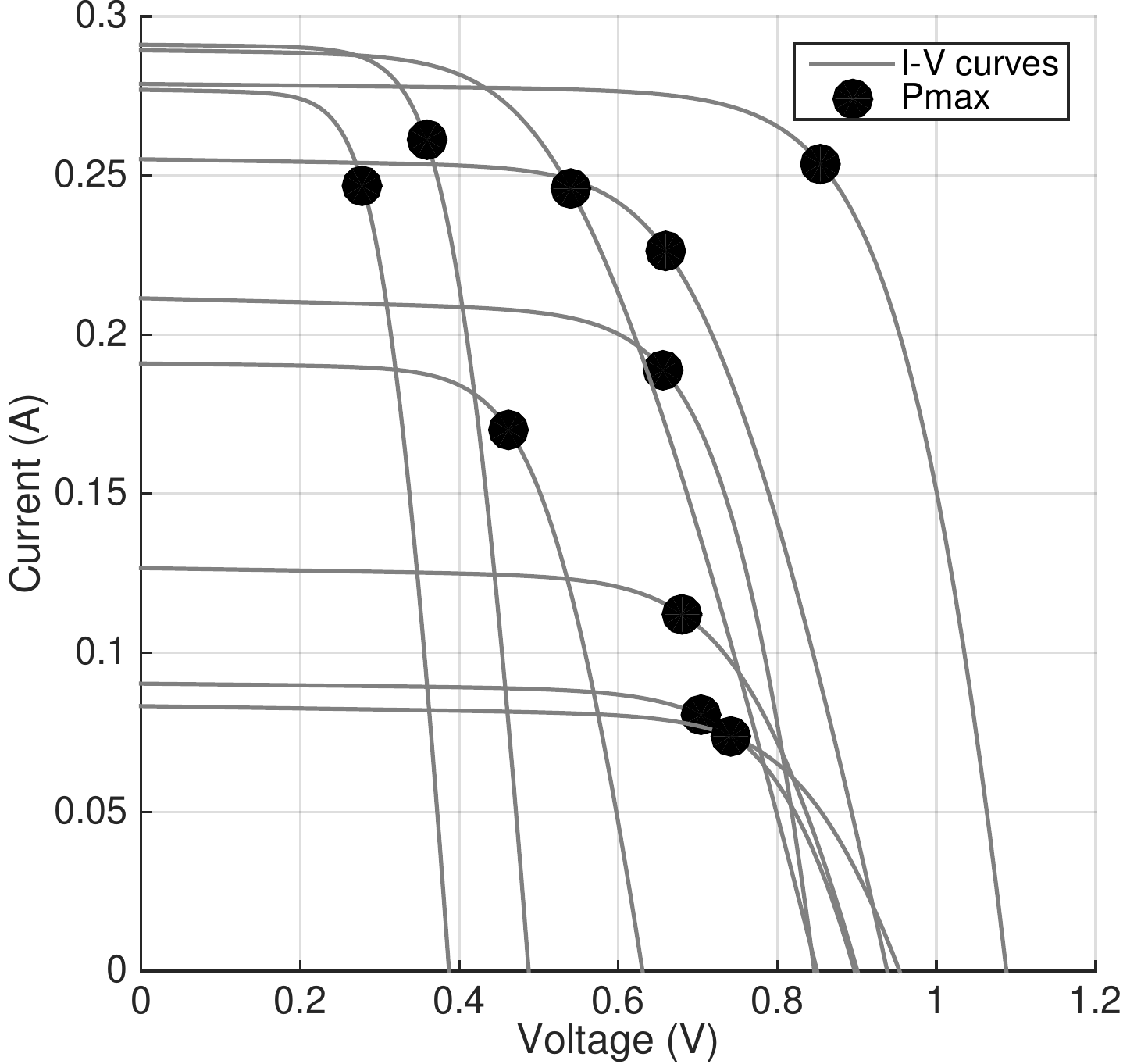}
\end{center}
\caption{Ten current-voltage curves with parameters chosen uniformly at random from the ranges in Table \ref{tab:domains}. The black circles plot the location of $\Pmax=I_{\text{max}}V_{\text{max}}$.}
\label{fig:ivcurves}
\end{figure}

Given a device and $I$-$V$ data, one can estimate precise values for the input parameters $(\ISC, \IS, n, \RS, \RP)$; see \cite{Hansen2013b,Farivar2010,Zaharatos2014}. With these values fixed, one can compute the maximum power of the device as 
\begin{equation}
\label{eq:power}
\Pmax \;=\; \maxi{I,V}\;I\, V,
\end{equation}
where current $I$ and voltage $V$ are constrained by \eqref{eq:model} and \eqref{eq:photocurrent}. Changing the inputs $(\ISC, \IS, n, \RS, \RP)$ changes the nonlinear constraints in the optimization \eqref{eq:power}, thus changing $\Pmax$. Therefore, we can write
\begin{equation}
\Pmax \;=\; \Pmax(\ISC, \IS, n, \RS, \RP). 
\end{equation}
We use MATLAB to solve the optimization \eqref{eq:power} subject to the constraints \eqref{eq:model} and \eqref{eq:photocurrent}. Figure \ref{fig:ivcurves} shows ten $I$-$V$ curves with the single-diode model parameters chosen uniformly at random from the ranges in Table \ref{tab:domains}. The black dots show the values of $I_{\text{max}}$ and $V_{\text{max}}$ that produce $\Pmax$. A wrapper function computes $\Pmax$ as a function of the input parameters. Requests for the MATLAB code may be sent to the third author Campanelli.  

\section{An active subspace in the single-diode model}
\label{sec:pvasm}

We apply the techniques from Section \ref{sec:asm} to search for an active subspace of the $\Pmax$ performance parameter as a function of the parameters defined in Table \ref{tab:domains}. Note that a different performance parameter (e.g., efficiency) would require an independent and identical analysis. The relationship between the respective active subspaces depends on the relationship between the performance parameters and their inputs. In some cases, the active subspaces might be similar; in other cases, they may differ substantially. 

The MATLAB scripts for generating the figures in this section can be found at \url{https://bitbucket.org/paulcon/active-subspaces-in-a-single-diode-model}. To generate the figures and results, we evaluated the model $(m+1)M=6000$ times, which took approximately three minutes on a MacBook Air. A more expensive model would permit fewer evaluations, and we would choose $M$ differently---or perhaps choose a different method.

\subsection{Normalizing the input parameters}

We first normalize the domain of $\Pmax$ to the hypercube; denote the normalized variables by $\vx=[x_1,x_2,x_3,x_4,x_5]^T$. The variable $\IS$ varies over several orders of magnitude, and preliminary tests show that $\Pmax$ changes rapidly near smaller values of $\IS$. To address this, we work with $\log(\IS)$. The range of $\log(\IS)$ is bounded below by -24.54 and above by -15.32. 

Let $\sL(\cdot)$ and $\sU(\cdot)$ return the upper and lower bounds, respectively, of the argument. We define the normalized input parameters as
\begin{equation}
\begin{aligned}
x_1 &= 2\,\left( \frac{\ISC-\sL(\ISC)}{\sU(\ISC)-\sL(\ISC)} \right) - 1\\
x_2 &= 2\,\left( \frac{\log(\IS)-\sL(\log(\IS))}{\sU(\log(\IS))-\sL(\log(\IS))} \right) - 1\\
x_3 &= 2\,\left( \frac{n-\sL(n)}{\sU(n)-\sL(n)} \right) - 1\\
x_4 &= 2\,\left( \frac{\RS-\sL(\RS)}{\sU(\RS)-\sL(\RS)} \right) - 1\\
x_5 &= 2\,\left( \frac{\RP-\sL(\RP)}{\sU(\RP)-\sL(\RP)} \right) - 1
\end{aligned}
\end{equation}
Then $\vx\in[-1,1]^5$ and has no units. Note that this map is invertible. In other words, given a value $\vx\in[-1,1]^5$, one can shift and scale the components of $\vx$---and transform the second component with the exponential---to produce a valid input for the single-diode model. We take the weight function $\rho$ in \eqref{eq:C} to be a constant $2^{-5}$ inside $[-1,1]$ and zero elsewhere. One can interpret this weight function as a uniform probability density on the space of input parameters. 

Proper scaling is important. The results of the active subspaces analysis are not scale invariant, so the scientist must carefully choose input parameter bounds that are appropriate for the application. Changing the bounds can change the results---sometimes dramatically. For example, the quantity of interest may be irregular and badly behaved over a wide range of inputs. But reducing the range of interest might focus on a small region in the parameter space where $f$ is smooth---even nearly linear. A linear function has a one-dimensional active subspace. Thus, loosely speaking, if the range of parameters is sufficiently small, then the active subspace will be nearly one-dimensional for continuous and differentiable quantities of interest. 

\subsection{Estimating the eigenvalues and eigenvectors}
We choose $M=1000$ points $\vx_i$ uniformly at random from $[-1,1]^5$, and for each $\vx_i$ we compute both $\Pmax$ and a first-order finite difference approximation of the gradient. We use a finite difference step size of $10^{-6}$ in the normalized domain, which was small enough to ensure that gradient approximation errors were substantially smaller than the variance due to random sampling.

\begin{figure}[h]
\begin{center}
\subfloat[Eigenvalues]{\label{fig:evals}
\includegraphics[width=0.45\textwidth]{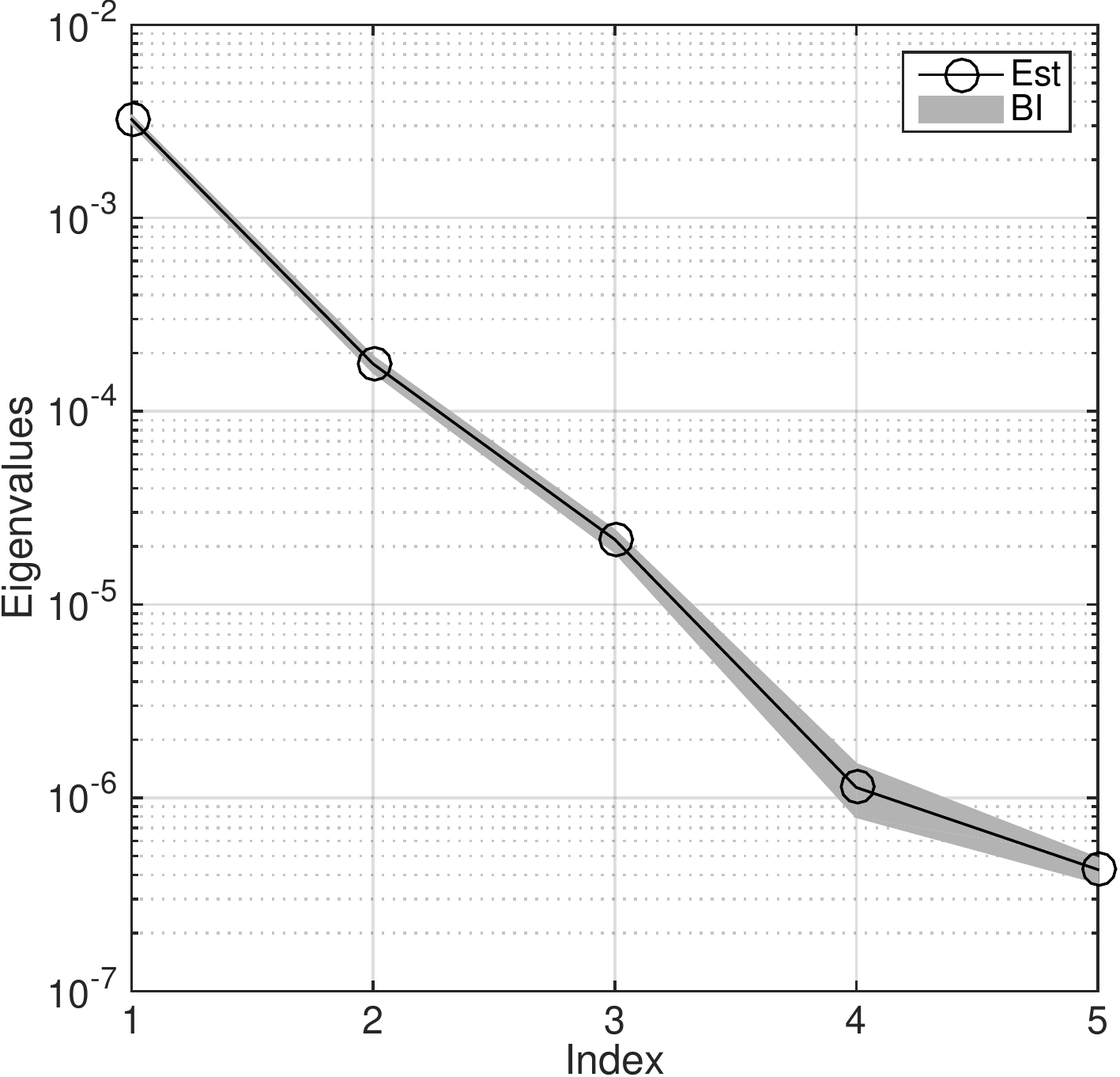}
}
\subfloat[Subspace error]{\label{fig:se}
\includegraphics[width=0.45\textwidth]{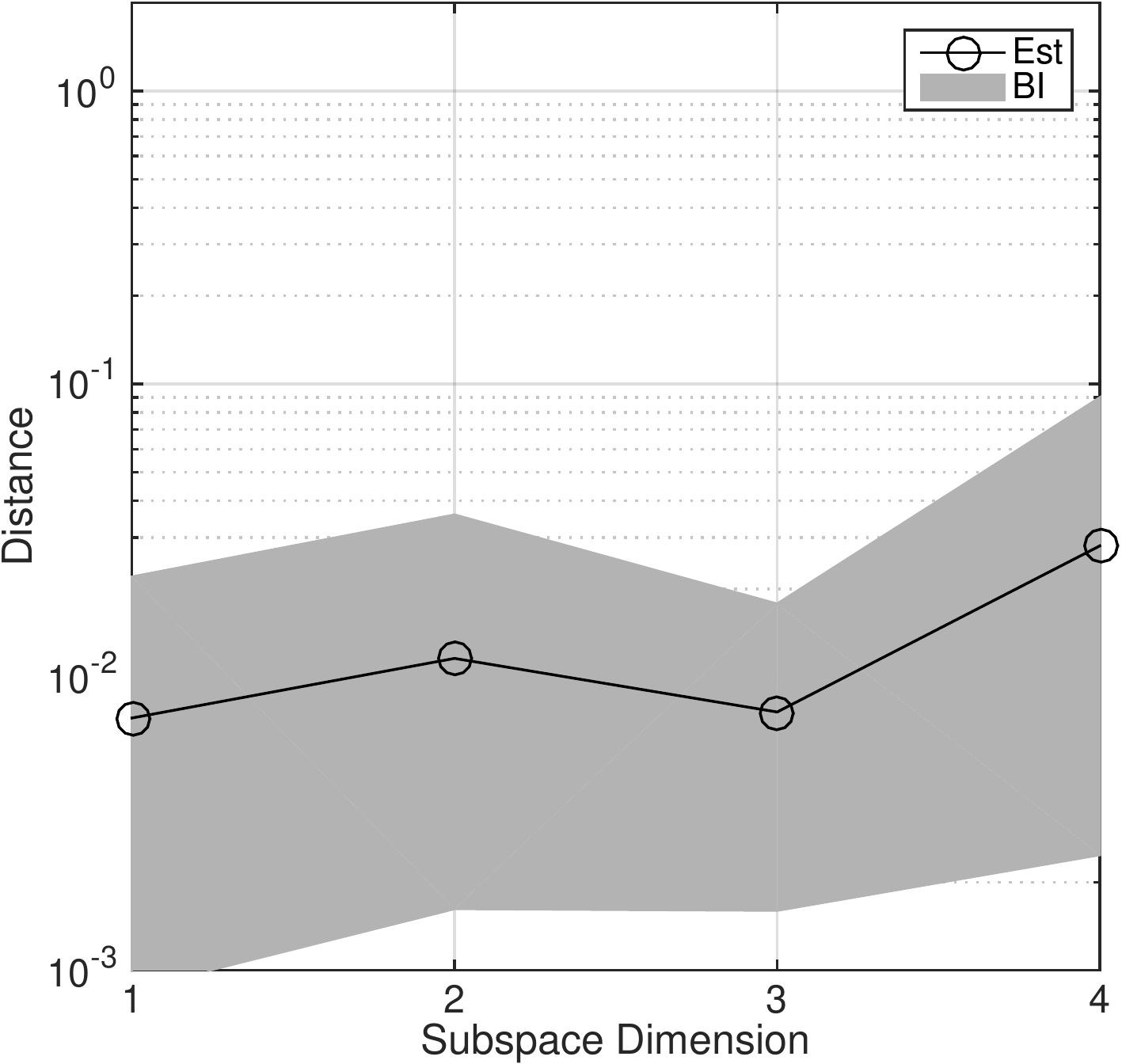}
}
\end{center}
\caption{The left figure shows the eigenvalues of the matrix $\hmC$ from \eqref{eq:Cmc} with $M=1000$ samples and their bootstrap ranges with $M'=1000$ boostrap replicates. The right figure shows the mean and ranges of the bootstrap estimates of the subspace error (see \eqref{eq:subspacerr}) with $M'=1000$ bootstrap replicates.}
 \label{fig:evals-and-se}
\end{figure}

To study variability in the computed components of $\hmC$ from \eqref{eq:Cmc}, we use the bootstrap described in Section \ref{sec:boot} with $M'=1000$ bootstrap replicates. The 99\% boostrap intervals for the eigenvalue estimates are shown in Figure \ref{fig:evals}. The tight bootstrap ranges suggest there is little eigenvalue variability in the bootstrapped data set. There is a gap of nearly an order of magnitude between both the first and second and the second and third eigenvalues. This suggests a very dominant one-dimensional active subspace. Figure \ref{fig:se} shows the estimated error in the approximated subspace using \eqref{eq:subspacerr}. Note that this error is actually the cosine of the principal angle between the subspaces, so it is bounded above by 1.

\begin{figure}[h]
\begin{center}
\subfloat[First eigenvector]{
\includegraphics[width=0.45\textwidth]{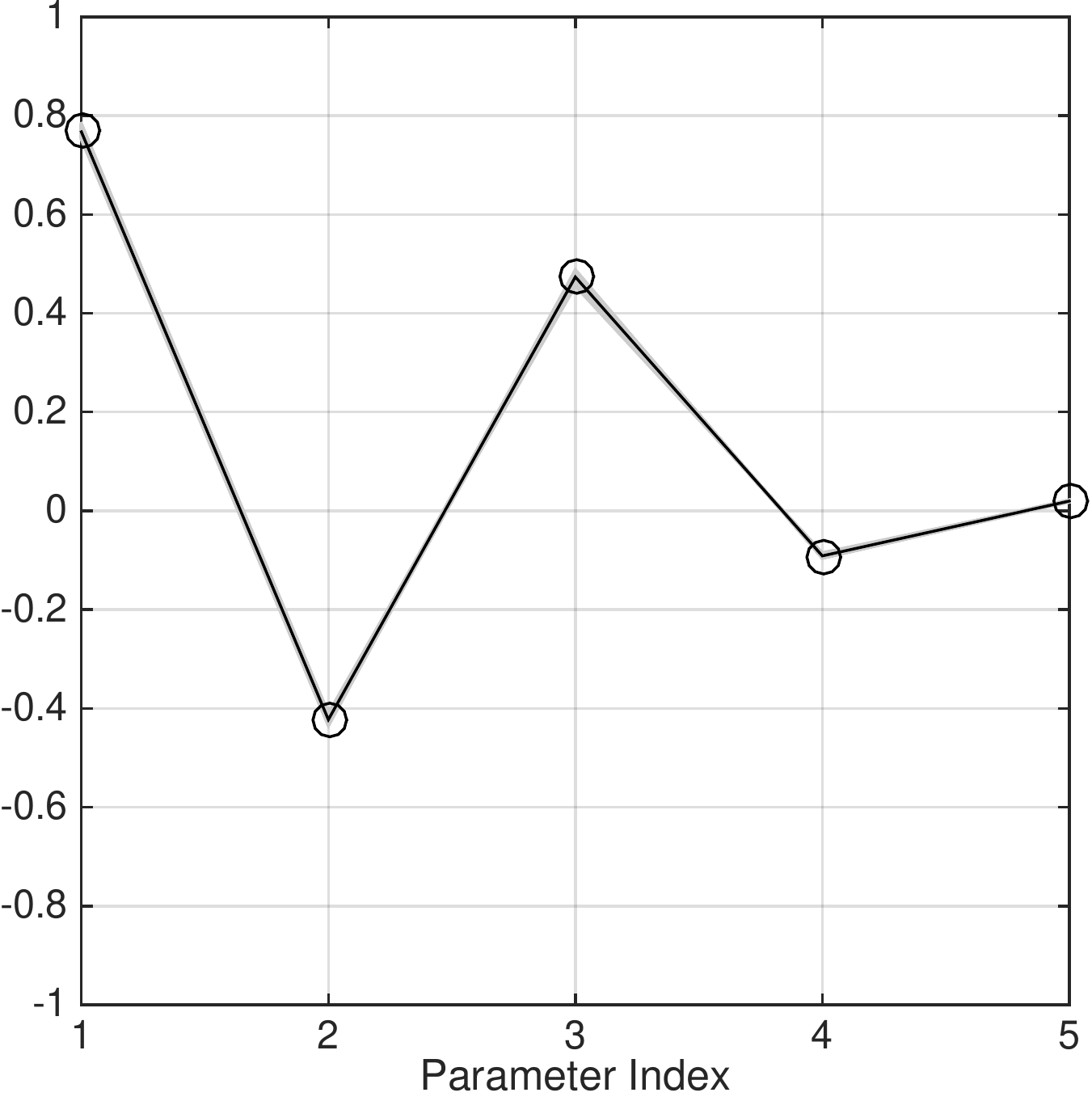}
}
\subfloat[Second eigenvector]{
\includegraphics[width=0.45\textwidth]{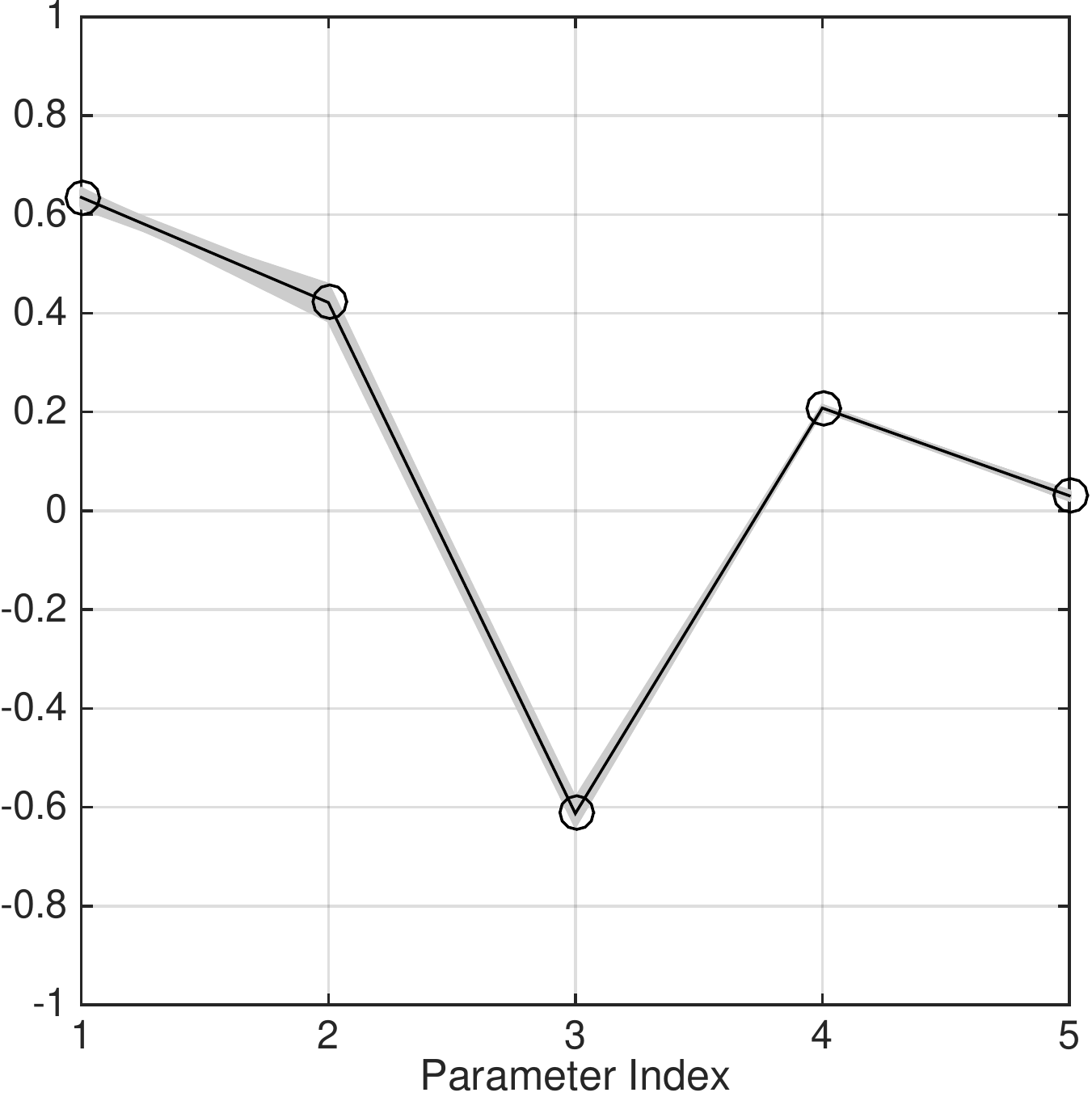}
}
\end{center}
\caption{These figures show the components of the first (left) and second (right) eigenvectors of $\hmC$ with $M=1000$ gradient samples. The grey shaded region superimposes the $M'=1000$ bootstrap replicates of the eigenvector components; the tight ranges suggests stability in the gradient samples. These values can be used as measures of sensitivity of $\Pmax$ with respect to the input variables in the model, whose names are given along the horizontal axis of the figures. See Table \ref{tab:domains} for a description of the variables.}
 \label{fig:evecs}
\end{figure}

Figure \ref{fig:evecs} displays the components of the first and second eigenvectors of $\hmC$ computed with $M=1000$ samples---along with the $M'=1000$ bootstrap replicates. They are normalized so that the first component is positive. The small ranges suggest that the one- and two-dimensional active subspaces are stable within the gradient samples. 

\subsection{Identifying the active subspace}
\begin{figure}[h]
\begin{center}
\subfloat[Projected with first eigenvector]{
\includegraphics[width=0.45\textwidth]{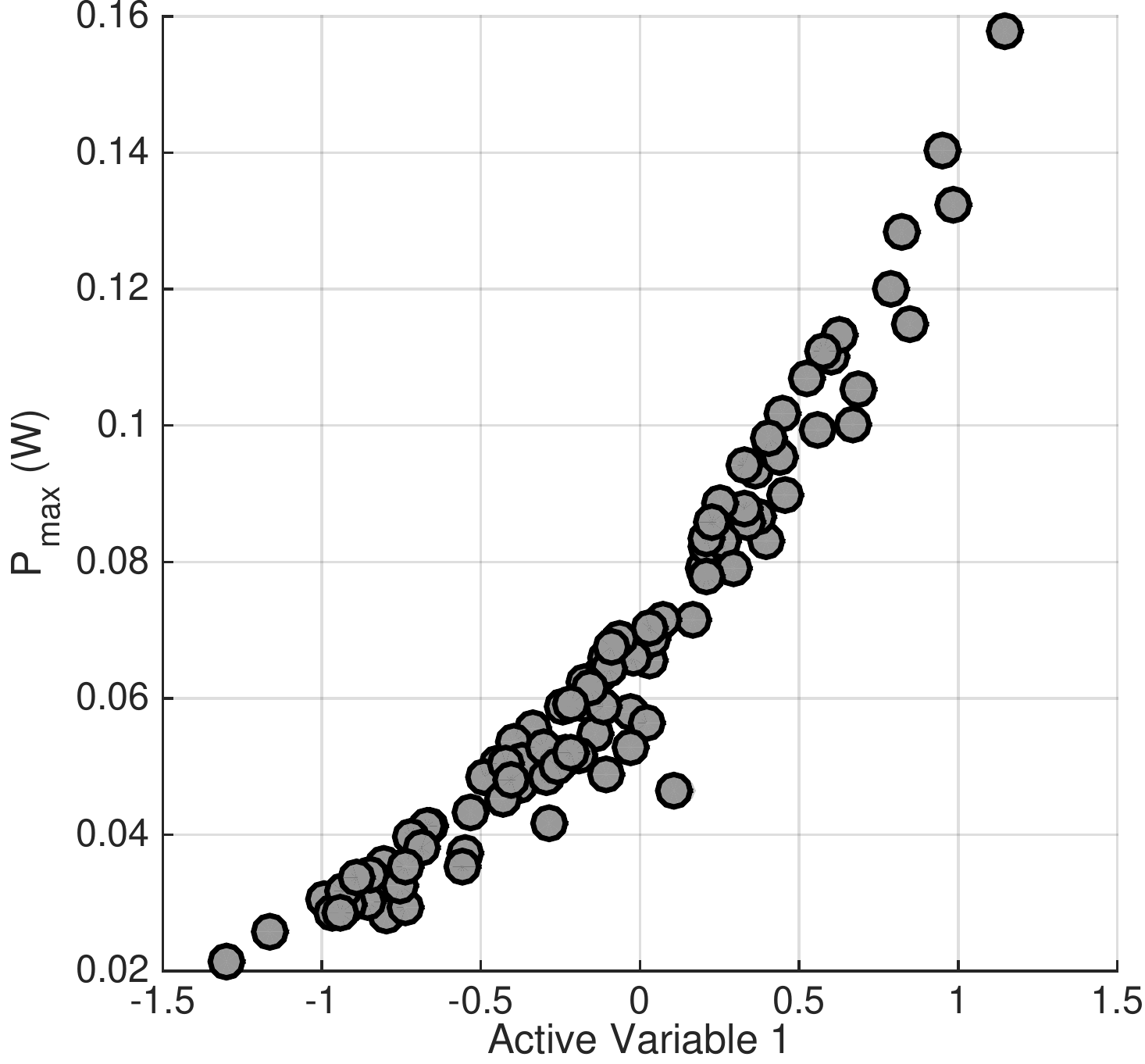}
\label{fig:oned}
}
\subfloat[Projected with first two eigenvectors]{
\includegraphics[width=0.45\textwidth]{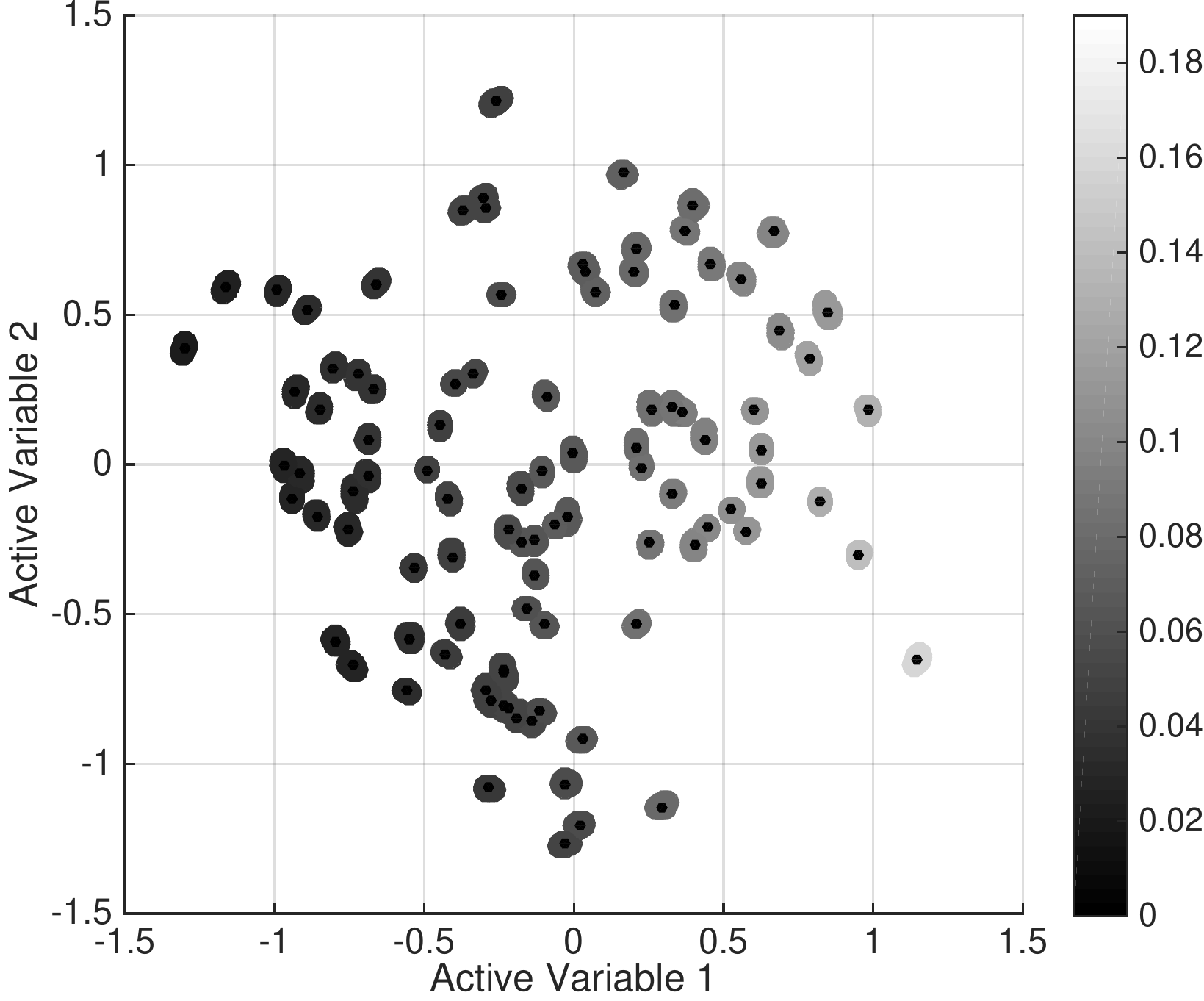}
\label{fig:twod}
}
\end{center}
\caption{One-dimensional (left) and two-dimensional sufficient summary plots of 100 realizations of $\Pmax$ using $\hvw_1$ (left) and both $\hvw_1$, $\hvw_2$ (right). The clusters of points on the right correspond to 20 bootstrap replicates of $\hvw_1$ and $\hvw_2$. The strong trend in \ref{fig:oned} is verified by looking at the variation along \emph{Active Variable 1} (the horizontal axis) in \ref{fig:twod}, and the combination of these provide evidence of a one-dimensional active subspace.}
\label{fig:proj1}
\end{figure}

The order-of-magnitude gap between the first and second eigenvalues suggests a dominant one-dimensional active subspace. Recall that gaps between eigenvalues are more important than the magnitudes of the eigenvalues. Figure \ref{fig:oned} plots 100 values of $\Pmax$ against a linear combination of the corresponding normalized inputs; the weights of the linear combination are the components of the first eigenvector $\hvw_1$ from \eqref{eq:Cmc}. Figure \ref{fig:oned} is a one-dimensional sufficient summary plot, as described in~\cite{cook2009regression}. The plot shows a potentially exploitable low-dimensional relationship between the input parameters and the quantity of interest, $\Pmax$. 

Figure \ref{fig:twod} shows the two-dimensional sufficient summary plot. The grayscale corresponds to the value of $\Pmax$. The black dot in the center of each cluster has a horizontal and vertical location determined by linear combinations of $\Pmax$'s corresponding input parameters; the weights of the linear combination are the eigenvectors $\hvw_1$ and $\hvw_2$, respectively. Each cluster corresponds to the value of $\Pmax$ with horizontal and vertical position determined by 20 bootstrap replicates of $\hvw_1$ and $\hvw_2$. The relatively small spread in the clusters indicates that the estimated subspace is stable. 

In Figure \ref{fig:twod}, note that $\Pmax$ increases primarily as a function of the first active variable. This provides more evidence of the dominance of the one-dimensional active subspace. 

\subsection{Eigenvectors for sensitivity analysis}

The eigenvector components connect the active subspace to the (normalized) variables in the model; Figure \ref{fig:evecs} includes the corresponding labels in the horizontal axis. The magnitudes of the eigenvector components can be used as measures of relative sensitivity for each of the parameters in the model. A large absolute value of the eigenvector component implies that this variable is important in defining the direction along which input perturbations produce the most change, on average, in $\Pmax$. The components corresponding the parallel resistance $\RP$ and series resistance $\RS$ are close to zero, which implies that (normalized) changes in $\RP$ and $\RS$ do not change $\Pmax$ as much as changes in the other parameters. 

We compare the components of $\hvw_1$, which defines the dominant one-dimensional active subspace, to the Sobol' first-order sensitivity indices and total sensitivity indices. We implemented the method to estimate these indices from Sudret~\cite{Sudret2008} using a tensor product Gaussian quadrature rule with 32768 points in 5 dimensions. Table \ref{tab:sens} displays the eigenvector components and the Sobol' indices. 

\begin{table}[ht]
\caption{The components of the eigenvector $\hvw_1$, the Sobol' first-order indices, and the Sobol' total sensitivity indices for the single-diode model parameters. The metrics provide same ranking of importance for the parameters.}
\begin{center}
\begin{tabular}{cccc}
Parameter & $\hvw_1$ & Sobol' first-order index & Sobol' total sens.~index\\
\hline
$\ISC$ & 0.77 & 0.56 & 0.61\\
$\IS$ & -0.42 & 0.17 & 0.19\\
$n$ & 0.47 & 0.21 & 0.25\\
$\RS$ & -0.09 & 0.01 & 0.01\\
$\RP$ & 0.02 & 0.00 & 0.00
\end{tabular}
\end{center}
\label{tab:sens}
\end{table}

There are several things to note in Table \ref{tab:sens}. First, the importance ranking from the Sobol' indices and the ranking from the eigenvector component magnitudes are the same. However, the signs of eigenvector components indicate whether $\Pmax$ will increase or decrease, on average, with changes in the corresponding parameter. (Eigenvectors are unique up to a sign, so these signs should be considered relative to one another.) There is no such interpretation in the Sobol' indices. In other words, the eigenvector components provide more information about the relationship between the corresponding parameters and $\Pmax$. 

The numeric values are difficult to compare, since they are normalized differently. The Sobol' indices are divided by the estimated total variance in $\Pmax$, whereas the eigenvector components are normalized to have Euclidean norm equal to 1. In this case, they are the same order of magnitude. 

We emphasize the different interpretations of these numbers. The eigenvector $\hvw_1$ is a single direction in the input space. Input perturbations along this direction change $\Pmax$ more, on average, than perturbations orthogonal to this direction. The Sobol' indices indicate the proportion of $\Pmax$'s total variance attributable to the factors of a variance-based decomposition of $\Pmax$. In terms of dimension reduction, one might (i) use the Sobol' indices to conclude that only three of the five parameters were important and (ii) approximate $\Pmax$ as a function of those three parameters. In contrast, the active subspace is used to approximate $\Pmax$ as with a univariate function of a linear combination of all five parameters. In other words, the variance-based approach reduces the dimension from 5 to 3, and the active subspace approach reduces the dimension from 5 to 1. 

There are other challenges facing the active subspace style of dimension  that do not affect the variance-based approach. The variance-based approach uses a subset of the model's parameters---as opposed to a subspace. Fixing a subset of parameters at nominal values and allowing the others to vary is relatively straighforward. Exploiting one important linear combination of the normalized parameters is not as straightforward. In the next section we discuss some possibilities for exploiting the active subspace. 

\subsection{Interpreting the sensitivity analyses}
We can interpret the consistent importance rankings in the sensitivity metrics in terms of the physics of the single-diode model. First, $\Pmax= I_{\text{max}}V_{\text{max}}$, where $(I_{\text{max}},V_{\text{max}})$ is the point on the $I$-$V$ curve that maximizes \eqref{eq:power}. $I_{\text{max}}$ typically scales linearly with $\ISC$. The voltage at $I=0$---denoted $V_{\text{oc}}$ for \emph{voltage at open circuit}---typically scales logarithmically with $\ISC$. $V_{\text{max}}$ behaves like $V_{\text{oc}}$. Thus, increases in $\ISC$ affect $I_{\text{max}}$ linearly and $V_{\text{max}}$ logarithmically, which changes $\Pmax$ as the product of a linear term and a logarithmic term. Since $\ISC$ affects both $I_{\text{max}}$ and $V_{\text{max}}$, it is natural for it to be the driving parameter in $\Pmax$. 

Recall that these sensitivity analyses are not scale invariant; the input parameter ranges and the weight function affect the results. The ranges for $n$ and $\IS$ are at about their maximum extent for a crystalline silicon PV cell, while the ranges for $\RS$ and $\RP$ are somewhat restricted. All of these parameters affect the shape and/or onset of the \emph{knee} of the $I$-$V$ curve, which determines the location of $(V_{\text{max}}, I_{\text{max}})$ and, thus, $\Pmax$. Here $\RS$ and $\RP$ are relatively unimportant because their ranges are relatively smaller than the ranges for $n$ and $\IS$.
 
\section{Using the active variables}
\label{sec:discussion}

The estimated eigenvalues in Figure \ref{fig:evals} and the plots in Figures \ref{fig:oned} and \ref{fig:twod} provide strong evidence for the presence of a dominant one-dimensional active subspace in the $\Pmax$ performance parameter computed from the single-diode model as a function of the model's five input parameters. The natural question is how one can exploit the dimension reduction afforded by the active subspace. There are five types of studies that benefit greatly from fewer input parameters. To keep the scope of this paper limited, we do not address the details of any of these studies, and we prefer to reserve them for future exploration. 

\noindent \textbf{Visualization.} The active subspace enables one to view the model output's dependence on its inputs with standard computer graphics tools when the active subspace is not more than two-dimensional. Sufficient summary plots, such as those in Figures \ref{fig:oned} and \ref{fig:twod} can provide insights to modelers seeking to improve their models. 

\noindent \textbf{Optimization.} Suppose one wanted to maximize $\Pmax$ over the input variables. The plot in Figure \ref{fig:proj1} shows a monotonic trend in the univariate function of the first active variable. For such functions, maximization is trivial; simply make the first active variable $\hvw_1^T\vx$ as large as possible subject to the constraints $\vx\in [-1,1]^m$. In general, a five-dimensional global optimization where the input/output relationship is not well-understood is a difficult problem. For the $\Pmax$ output from the single-diode model, the active subspace provides a way to discover the location of the global optimum with ease. 

\noindent \textbf{Response surfaces.} If the single-diode model were expensive to evaluate, then one may wish to construct a response surface that approximates the map from inputs to outputs. A one-dimensional active subspace allows one to build a response surface on only the active variable instead of the five model input parameters. Constructing a response surface in one variable is certainly preferred to constructing one in five variables. Real-time control systems such as maximum power point trackers might benefit from the characterization of a simple, low-dimensional response surface. In addition, the visualization tools give one confidence that the response is sufficiently smooth with respect to the active variable to permit an accurate response surface. 

\noindent \textbf{Averages.} One may wish to compute an average of $\Pmax$ over all five input variables. Since the average of $\Pmax$ over its inputs is the same as the average of the conditional expectation given the active variable, we need only to approximate a marginal density of the active variable to compute the average. We are working on the specifics of this idea, but the goal would be to transform a five-dimensional integral into a one-dimensional integral requiring many fewer evaluations of $\Pmax$.

\noindent \textbf{Design.} Suppose a modeler wanted to design a solar cell with $\Pmax$ in a specified range and sought an appropriate range of values for the model inputs. A one-dimensional active subspace makes this query much easier---especially if $\Pmax$ is monotonic with respect to the active variable as in Figure \ref{fig:twod}. 

\section{Conclusions}
\label{sec:conclusion}

We have discussed methods for discovering the directions in a model's input space that change the model's prediction the most, on average. The active subspace is the span of these directions. For global sensitivity analysis, the active subspace offers an alternative to variance decomposition techniques such as Sobol' indices. We have applied this procedure to a single-diode solar cell model with key performance parameter $\Pmax$, and we discovered a dominant one-dimensional active subspace. We offered several possible ways to exploit the knowledge of this low-dimensional parametric dependence to gain greater insight into the model.

\pagebreak

\bibliographystyle{natbib}
\bibliography{pv-asm-coda}

\end{document}